\input amstex
\documentstyle{amsppt}
\document
\magnification=1200
\NoBlackBoxes
\nologo
\pageheight{18cm}


\bigskip

\centerline{\bf GEORG CANTOR AND HIS HERITAGE\footnotemark1}
\footnotetext{Talk at the meeting of the German Mathematical Society
and the Cantor Medal award ceremony.}

\medskip

\centerline{\bf Yu.~I.~Manin}

\bigskip

{\it God is no geometer, rather an unpredictable poet.}

{\it (Geometers can be unpredictable poets,
so there could be room for compromise.)}

\smallskip

\hfill{V.~Tasi\'c [T], on the XIXth century romanticism}

\bigskip

\centerline{\bf Introduction}

\medskip

Georg Cantor's grand meta-narrative, Set Theory,
created by him almost singlehandedly in the span
of about fifteen years, resembles a piece 
of high art more than a scientific theory.

\smallskip

Using a slightly modernized language,
basic results of set theory can be stated in a few lines.

\smallskip

Consider the category of all sets with arbitrary
maps as morphisms. Isomorphism classes
of sets are called {\it cardinals}. Cardinals are well--ordered
by the sub--object relation, and the
cardinal of the set of all subsets of $U$
is strictly larger than that of $U$ (this
is of course proved by the famous diagonal argument).

\smallskip

This motivates introduction of another category,
that of well ordered sets and monotone maps as morphisms.
Isomorphism classes of these are called
{\it ordinals}. They are well--ordered as well.
The Continuum Hypothesis is a guess about the
order structure of the initial segment of cardinals.

\smallskip

Thus, exquisite minimalism of expressive means is used by Cantor to 
achieve a sublime goal: understanding  infinity,
or rather infinity of infinities. A built--in 
self--referentiality and the forceful extension
of the domain of mathematical intuition
(principles for building up new sets) add
to this impression of combined artistic violence and self--restraint.

\smallskip

Cantor himself would have furiously opposed this view.
For him, the discovery of the hierarchy of infinities
was a revelation of God--inspired Truth.

\smallskip

But mathematics of the XXth century  reacted to Cantor's oeuvre
in many ways that can be better understood in the
general background of various currents of
contemporary science, philosophical thought, and
art.

\smallskip

Somewhat provocatively, one can render 
one of Cantor's principal  insights as follows:

\medskip
\centerline{\it $2^x$ is considerably larger than $x$.}
\medskip

Here $x$ can be understood as an integer, an arbitrary ordinal,
or a set; in the latter case $2^x$ denotes the set of all subsets
of $x$.
Deep mathematics starts when we
try to make this statement more precise
and to see {\it how much larger} $2^x$ is. 

\smallskip

If $x$ is the first infinite ordinal,
this is the Continuum Problem. 

\smallskip

I will argue
that properly stated for finite $x$, this question becomes
closely related to a universal $NP$ problem.

\smallskip

I will then discuss assorted topics related to the role of set theory
in contemporary mathematics and the reception of Cantor's ideas.

\bigskip

\centerline{\bf Axiom of Choice and $P/NP$--problem}

\centerline{\it or}

\centerline{\bf finite as poor man's transfinite\footnotemark2}
\footnotetext{I am not hinting at the Clay  USD $10^6$ Prize
for a solution of the $P/NP$--problem.}

\medskip

In 1900, at his talk at the second ICM
in Paris, Hilbert put the Continuum Hypothesis
at the head of his list of 23 outstanding
mathematical problems. This was one of the highlights
of Cantor's scientific life. Cantor invested much effort consolidating the
German and international mathematical community
into a coherent body capable of counterbalancing
a group of influential professors
opposing set theory.   
 
\smallskip

The opposition to his theory of infinity, however, continued and
was very disturbing to him, because the validity of Cantor's new mathematics was
questioned. 

\smallskip

In 1904, at the next International Congress,
K\"onig presented a talk in which he purported to show that 
continuum could not be well--ordered,
and therefore the Continuum Hypothesis was meaningless.

\smallskip

Dauben writes ([D], p. 283): ``The dramatic events of K\"onig's paper read during the
Third International Congress for Mathematicians in Heidelberg greatly upset
him [Cantor]. He was there with his two daughters, Else and Anna--Marie, and was outraged
at the humiliation he
felt he had been made to suffer.''

\smallskip

It turned out that K\"onig's paper contained a mistake; 
Zermelo soon afterwards produced a proof that
any set could be well ordered using his then brand new Axiom of Choice.
This axiom essentially postulates that, starting with a set
$U$, one can form a new set, whose elements
are pairs $(V,v)$ where $V$ runs over all non--empty subsets
of $U$, and $v$ is an element of $V$.

\smallskip

A hundred years later, the mathematical community 
did not come up with a compendium of new problems for the coming
century similar to the Hilbert's list. Perhaps, the general vision
of mathematics changed --- already in Hilbert's list
a considerable number of items could be better described as research programs
rather than well--defined problems, and this seems to be
a more realistic way of perceiving our work in progress.

\smallskip

Still, a few clearly stated important questions remain unanswered, and recently
seven such questions were singled out and endowed with a price tag.
Below I will invoke one of these questions, the $P/NP$ problem 
and look at it as  a finitary travesty
of Zermelo's Axiom of Choice.

\smallskip

Let $U_m=\bold{Z}_2^m$ be the set of
$m$--bit sequences. A convenient way to encode its subsets
is via Boolean polynomials. Using the standard -- and more
general -- language of commutative algebra,
we can identify
each subset of $U_m$ 
with the $0$--level of a unique function $f\in B_m$
where we define the algebra of Boolean polynomials as
$$
B_m:=\bold{Z}_2[x_1,\dots ,x_m]/(x_1^2+x_1, \dots ,x_m^2+x_m).
$$
Hence Zermelo's problem -- {\it choose an element in each non--empty subset
of $U$} -- translates into: {\it for each Boolean polynomial,
find a point at which the polynomial equals 0, or prove that
the polynomial is identically 1}. Moreover, we want to solve
this problem in time bounded by a polynomial
of the bit size of the code of $f$.

\smallskip

This leads to a universal
(``maximally difficult'') $NP$--problem if one writes 
Boolean polynomials in the following
version of disjunctive normal form. Code of such a
form is a family
$$
u =\{m; (S_1,T_1),\dots ,(S_N, T_N)\}, \ m\in \bold{N}; \, S_i,T_i\subset
\{1,\dots ,m\}.
$$
The bit size of  $u$ is $mN$, and the respective Boolean polynomial is
$$
f_u:=1+ \prod_{i=1}^N\left(1+\prod_{k\in S_i}(1+x_k)\prod_{j\in T_i}x_j\right)
$$
This encoding provides for a fast check of the inclusion relation
for the elements of the respective level set.  The price
is that the uniqueness
of the representation of $f$ is lost, and moreover,
the identity $f_u=f_v (?)$ becomes a computationally hard problem.

\smallskip

 In particular,
even the following weakening of the finite Zermelo problem
becomes $NP$--complete and hence currently intractable:
{\it check whether a Boolean polynomial given in  disjunctive normal form is
non--constant.}

\smallskip 

Zermelo's Axiom of Choice aroused a lively international discussion 
published in the first issue of Mathematische Annalen of 1905.
A considerable part of it was focused on the psychology of mathematical 
imagination and on the reliability of its fruits.
Baffling questions of the type: ``How can we be sure that during the course
of a proof we keep thinking about the same set?'' kept
emerging. If we imagine that at least a part of computations
that our brains perform can be adequately modeled by finite automata,
then quantitative estimates of the required resources as
provided by the theory of polynomial time computability
might eventually be of use in neuroscience and by implication in psychology.

\smallskip

A recent paper in ``Science''  thus summarizes some
experimental results throwing light on the nature
of mental representation of mathematical objects and physiological roots
of divergences between, say, intuitionists and formalists:

\smallskip

``[...] our results provide grounds for reconciling the divergent
introspection of mathematicians by showing that even within the small
domain of elementary arithmetic, multiple mental
representations are used for different tasks. Exact arithmetic puts emphasis
on language--specific representations and relies on a left inferior frontal
circuit also used for generating associations between words.
Symbolic arithmetic is a cultural invention specific to humans, and its development
depended on the progressive improvement of number
notation systems. [...]

Approximate arithmetic, in contrast, shows no dependence
on language and relies primarily on a quantity representation implemented
in visuo--spatial networks of the left and right
parietal lobes.'' ([DeSPST], p.~973).

\smallskip

In the next section, I will discuss an approach
to the Continuum Conjecture which is clearly inspired
by the domination of the visuo--spatial networks,
and conjecturally better understood in terms
of probabilistic models than logic or Boolean
automata.

\smallskip

{\it Appendix: some definitions.} For completeness, I will remind 
the reader of the basic definitions
related to the $P/NP$ problem. 
Start with an infinite constructive world $U$ in the sense
of [Man], e.g. natural numbers $\bold{N}$.
A subset $E\subset U$ {\it belongs to the class $P$}
if it is decidable  and the values of its characteristic function $\chi_E$ are computable in
polynomial time on all arguments $x\in E.$

\smallskip

Furthermore, $E\subset U$ {\it belongs to the class $NP$}
if it is a polynomially truncated projection
of some $E'\subset U\times U$ belonging to $P$:
for some polynomial $G$,
$$
u\in E\ \Leftrightarrow\ \exists\,(u,v)\in E'\ \roman{with}\ |v|\le G(|u|)
$$
where $|v|$ is the bit size of $v$. In particular, $P\subset NP$.

\smallskip

Intuitively, $E\in NP$ means that for each $u\in E$ {\it there exists
a polynomially bounded proof} of this inclusion
(namely, the calculation of $\chi_{E'} (u,v)$ for
an appropriate $v$). However, to find such a proof (i.~e. $v$)
via the naive search among all $v$ with $|v|\le G(|u|)$ 
can take exponential time.

\smallskip

The set $E\subset U$ is called {\it NP--complete} if, for any other
set $D\subset V, D\in NP,$ there exists a  polynomial time
computable function $f:\,V\to U$ such that $D=f^{-1}(E)$,
that is, $\chi_D(v)=\chi_E(f(v)).$ 

\smallskip

The  encoding of Boolean polynomials used above is
explained and motivated by the proof of the
$NP$--completeness: see e.~g. [GaJ], sec. 2.6.

\bigskip


\centerline{\bf The Continuum Hypothesis and random variables}

\medskip

Mumford in [Mum], p. 208, recalls an argument of Ch.~Freiling ([F])
purporting to show that Continuum Hypothesis is ``obviously'' false
by considering the following situation: 

\smallskip

``Two dart players independently throw darts
at a dartboard. If the continuum hypothesis is true,
the points $P$ on the surface of a dartboard can be well ordered so that
for every $P$, the set of $Q$ such that $Q<P$, call it $S_Q$, is countable.
Let players 1 and 2 hit the dart board
at points $P_1$ and $P_2$. Either
$P_1<P_2$ or $P_2<P_1$. Assume the first holds. Then $P_1$ belongs to a countable subset
$S_{P_2}$ of the points on on the dartboard. As the two
throws were independent, we may treat throw 2 as taking place first,
then throw 1. After throw 2, this countable set $S_{P_2}$
has been fixed. But every countable
set is measurable and has measure 0. The same argument shows that
the probability of $P_2$ landing
on $S_{P_1}$ is 0. Thus almost surely
neither happened and this contradicts
the assumption that the dartboard is the first
uncountable cardinal! [...]

\smallskip

I believe [...] his `proof' shows that if we make random variables one of the
basic elements of mathematics, it follows that
the C.H. is false and we will get rid of
one of the meaningless conundrums of set theory.''\footnotemark3
\footnotetext{Mumford's paper is significantly
called {\it ``The dawning of the age of stochasticity.''} David assured me
that he had no intention to refer to Giambattista Vico's
theory of historical cycles which in the rendering of  H.~Bloom ([B])
sounds thus: `` Giambattista Vico, in his {\it New Science}, posited a cycle
of three phases --- Theocratic, Aristocratic, Democratic ---
followed by a chaos of which a New Theocratic Age would at last emerge.''}

\smallskip

Freiling's work was actually preceded
by that of Scott and Solovay, who recast
in terms of ``logically random sets''
P.~Cohen's forcing method for proving the consistency
of the negated CH with Zermelo--Frenkel axioms.
Their work has shown that one can indeed put random
variables in the list of basic notions
and use them in a highly non--trivial way.

\smallskip
 
P.~Cohen himself ended his book by
suggesting that the view that CH is ``obviously false''
may become universally accepted.

\smallskip

However, whereas the Scott--Solovay reasoning
proves a precise theorem about the formal language
of the Set Theory, Freiling's argument
appeals directly to our physical intuition, and 
is best classified
as a {\it thought experiment.} It is similar in nature
to some classical thought experiments in physics,
deducing e.g. various dynamic consequences from
the impossibility of {\it perpetuum mobile}.

\smallskip
 
The idea of thought experiment, as opposed to that of
logical deduction, can be generally considered as 
a right--brain equivalent of the left brain elementary
logical operations. A similar role is played
by good metaphors. When we are comparing the respective 
capabilities of two brains, we are struck by what I called elsewhere
`the inborn weakness of metaphors': they resist  becoming building blocks
of a  system. One can only more or less artfully put one metaphor upon another, and
the building will stand or crumble upon its own weight
independently of its truth or otherwise. 

\smallskip

Physics disciplines
thought experiments as poetry disciplines metaphors,
but only logic has an inner discipline.

\smallskip

Successful thought experiments produce mathematical truths which,
after being accepted,
solidify into axioms, and the latter start working
on the treadmill of logical deductions.

\bigskip

\centerline{\bf Foundations  and Physics}

\medskip

I will start this section with a brief discussion of
the impact of Set Theory on the foundations of mathematics. 
I will understand ``foundations'' neither as the para--philosophical 
preoccupation with nature, accessibility, and reliability 
of mathematical truth, nor as a set of normative prescriptions 
like those advocated by finitists or formalists.

\smallskip

I will use this word in a loose sense as
a general term for the  historically
variable conglomerate of rules and principles
used to organize the already existing and always being created anew
body of mathematical knowledge of the relevant epoch.
At times, it becomes codified in the form of 
an authoritative mathematical text as exemplified by
Euclid's Elements. In another epoch, it is better expressed
by nervous self--questioning about the meaning of
infinitesimals or the precise relationship between real numbers
and points of the Euclidean line, or
else, the nature of algorithms. In all cases, foundations in this wide sense
is something which is of relevance to a working
mathematician, which refers to some basic principles of
his/her trade, but which does not constitute the essence of his/her work. 

\smallskip

In the XXth century, all of the main foundational trends referred to 
Cantor's language and intuition of sets. 

\smallskip

The well developed  project of Bourbaki gave a polished form
to the notion that {\it any} mathematical object $\Cal{X}$
(group, topological space, integral, formal language ...)
could be thought of as a set $X$ with an additional structure $x$.
This notion had emerged in many 
specialized research projects, from Hilbert's
{\it Grundlagen der Geometrie} to
Kolmogorov's identification of probability
theory with measure theory.

\smallskip

The additional structure $x$ in $\Cal{X}=(X,x)$
is an element of another set $Y$ belonging to the {\it \'echelle} constructed
from $X$ by standard operations and satisfying conditions (axioms)
which are also formulated entirely in terms of set theory.
Moreover, the nature of elements of $X$ is inessential:
a bijection $X\to X'$ mapping $x$ to $x'$ produces an isomorphic object
$\Cal{X}'=(X',x')$. This idea played a powerful unifying and clarifying role in mathematics
and led to spectacular developments far outside the Bourbaki group.
Insofar as it is accepted in thousands of research papers,
one can simply say that the language of mathematics is
the language of set theory.

\smallskip

Since the latter is so easily formalized, this allowed logicists
to defend the position that their normative principles
should be applied to all of mathematics and to overstate
the role of ``paradoxes of infinity'' and G\"odel's incompleteness
results.

\smallskip
 
However, this fact also made possible such self--reflexive
acts as inclusion of metamathematics into mathematics, in the form of
model theory. Model theory studies special algebraic structures -- formal languages --
considered in turn as mathematical objects (structured sets
with composition laws, marked elements etc.),
and their interpretations
in sets. Baffling discoveries such as G\"odel's incompleteness of arithmetics
lose some of their mystery once one comes to understand their
content as a statement that a certain algebraic structure
simply is not finitely generated with respect to the allowed composition laws.

\smallskip

When at the next stage of this historical development, sets gave way to
categories,  this was at first only a shift of stress upon morphisms
(in particular, isomorphisms) of structures, rather than on structures
themselves. And, after all, a (small) category could itself
be considered as a set with structure. However, primarily thanks to the work
of Grothendieck and his school on the foundations of algebraic geometry, 
categories moved to the foreground. Here is an incomplete list
of changes in our understanding of mathematical objects brought about
by the language of categories. Let us recall that generally
objects of a category $C$ are not sets themselves; their nature
is not specified; only
morphisms $Hom_C(X,Y)$ are sets.

\smallskip

A. An object $X$ of the category $C$ can be identified
with the functor it represents: $Y\mapsto Hom_C(Y,X)$.
Thus, if $C$ is small, initially structureless $X$ becomes a structured set. This external, ``sociological''
characterization of a mathematical object defining it
through its interaction with all
objects of the same category rather than in terms of its
intrinsic structure, proved to be extremely useful
in all problems involving e.~g. moduli spaces in algebraic geometry.

\smallskip

B. Since two mathematical objects, if they are isomorphic,
have exactly the same properties, it does not matter
how many pairwise isomorphic objects are contained
in a given category $C$. Informally, if
$C$ and $D$ have ``the same'' classes of isomorphic
objects and morphisms between their representatives,
they should be considered as equivalent. For example,
the category of ``all'' finite sets is equivalent to any
category of finite sets in which there is exactly one
set of each cardinality 0,1,2,3, $\dots$\, .

\smallskip

This ``openness'' of a category considered up to equivalence
is an essential trait, for example, in the abstract computability
theory. Church's thesis can be best understood
as a postulate that there is an open category
of ``constructive worlds'' --- finite or countable structured sets 
and computable morphisms between them ---
such that any infinite object in it is isomorphic to the world of
natural numbers, and morphisms  correspond to recursive functions
(cf. [Man] for more details).  There are many more
interesting infinite constructive worlds determined
by widely diverging internal structures: 
words in a given alphabet, finite graphs, 
Turing machines, etc. However, they are all isomorphic
to $\bold{N}$ due to the existence of computable numerations.

\smallskip

C.  The previous remark also places limits on the naive view that
categories ``are'' special structured sets. In fact,
if it is natural to identify categories related by an equivalence
(not necessarily bijective on objects) rather than isomorphism, then this view becomes
utterly misleading.

\smallskip

More precisely, what happens is the slow emergence of the following
hierarchical picture. Categories themselves form objects
of a larger category $Cat$ morphisms in which are functors,
or ``natural constructions'' like a (co)homology theory of topological
spaces. However, functors do not form simply a set
or a class: they also form objects of a category.  Axiomatizing
this situation we get a notion of {\it 2--category} whose
prototype is $Cat$. Treating 2--categories in the same way, we get 3--categories etc.

\smallskip

The following view of mathematical objects is encoded in this hierarchy:
there is no equality of mathematical objects, only
equivalences. And since an equivalence is also a mathematical object,
there is no equality between them, only the next order equivalence etc., {\it ad infinitim}.

\smallskip

This vision, due initially to Grothendieck, extends the boundaries
of classical mathematics, especially algebraic geometry, 
and exactly in those developments where it interacts
with modern theoretical physics.

\smallskip

With the advent of categories,
the mathematical community was  cured of its fear of
classes (as opposed to sets) and generally ``very large''
collections of objects. 

\smallskip

In the same vein, it turned out that there are meaningful ways of thinking
about ``all'' objects of a given kind, and to use self--reference creatively
instead of banning it completely. This is a development
of the old distinction between sets and classes,
admitting that at each stage we get a structure similar to but
not identical with the ones we studied at the previous stage.

\smallskip

In my view, Cantor's prophetic vision was
enriched and not shattered by these new developments.

\smallskip

What made it recede to the background, together with
preoccupations with paradoxes of infinity
and intuitionistic neuroses,
was a renewed interaction with physics and the  
transfiguration of formal logic into computer science.

\smallskip

The birth of quantum physics radically changed
our notions about relationships between reality, its theoretical descriptions,
and our perceptions. It made clear that 
Cantor's famous definition of sets ([C]) represented only
a distilled classical mental view of the material world as consisting
of pairwise distinct things residing in space:

\smallskip

{\it ``Unter einer  `Menge'
verstehen wir jede Zusammenfassung $M$ von bestimmten
wohlunterschiedenen Objekten $m$ unserer Anschauung
oder unseres Denkens (welche die `Elemente' von $M$ genannt werden)
zu einem Ganzen.''}

\smallskip

``By a `set' we mean any collection $M$ into a whole of definite,
distinct objects $m$ (called the `elements' of $M$)
of our perception  or our thought.''

\smallskip

Once this view was shown to be only an approximation
to the incomparably more sophisticated quantum description,
sets lost their direct roots in reality.
In fact, the structured sets of modern mathematics used most
effectively in modern physics are not sets of things,
but rather sets of {\it possibilities}. For example, the phase space
of a classical mechanical system consists
of  pairs {\it (coordinate, momentum)}
describing all possible states of the system, whereas after quantization
it is replaced by the space of complex probability amplitudes:  
the Hilbert space of $L_2$--functions of the
coordinates or something along these lines. The amplitudes are
all possible quantum superpositions of all
possible classical states. It is a far cry from
a set of things.

\smallskip

Moreover, requirements of quantum physics
considerably heightened the degree of tolerance  
of mathematicians to the imprecise
but highly stimulating theoretical discourse
conducted by physicists.
This led, in particular, to the emergence
of Feynman's path integral as one of the most active areas of
research in topology and algebraic geometry, even though
the mathematical status of path integral is in no better shape
than that of Riemann integral before Kepler's ``Stereometry of Wine Barrels''.

\smallskip

Computer science added a much needed touch of practical relevance
to the essentially hygienic prescriptions of formal logic.
The introduction of the notion of ``success with high probability''
into the study of algorithmic solvability
helped to further demolish mental barriers which
fenced off foundations of mathematics from mathematics itself.  
 
\smallskip

{\it Appendix: Cantor and physics.} It would be interesting
to study Cantor's natural philosophy
in more detail. According to [D], he directly referred
to possible physical applications of his theory
several times. 

\smallskip

For example, he proved that 
 that if one deletes from a domain in
$\bold{R}^n$ any dense countable subset (e.g. all algebraic
points), then any two points of the complement can be
connected by a continuous 
curve. His interpretation: continuous motion
is possible even in discontinuous spaces, so ``our'' space 
might be discontinuous itself, because the idea
of its continuity is based upon observations
of continuous motion. Thus a revised mechanics should be considered.

\smallskip

At a meeting of GDNA in Freiburg, 1883, Cantor said: ``One of the
most important problems of set theory [...] consists
of the challenge to determine the various valences or powers of sets present
in all of Nature in so far as we can know them'' ([D], p. 291).

\smallskip

Seemingly, Cantor wanted atoms (monads) to be actual points, extensionless,
and in infinite number in Nature. ``Corporeal monads''
(massive particles? Yu.M.) should exist in countable
quantity. ``Aetherial monads'' (massless quanta? Yu.M.)
should have had cardinality aleph one.

\bigskip

\centerline{\bf Coda: Mathematics and postmodern condition}

\medskip

Already during Cantor's life time, the reception of his ideas 
was more like that of
new trends in the art, such as impressionism or atonality, 
than that of new scientific theories.
It was highly emotionally charged and ranged from
total dismissal (Kronecker's ``corrupter of youth'')
to highest praise (Hilbert's defense of ``Cantor's Paradise'').
(Notice however the commonly overlooked nuances
of both statements which subtly undermine their ardor:  
Kronecker implicitly likens  Cantor to Socrates, whereas Hilbert with faint mockery
hints at Cantor's conviction that  Set Theory
is inspired by God.)

\smallskip

If one accepts the view that Bourbaki's
vast construction was the direct descendant
of Cantor's work, it comes as no surprise
that it shared the same fate: see [Mas]. Especially
vehement was reaction against ``new maths'':
an attempt to reform the mathematical
education by stressing precise definitions, logic,
and set theoretic language rather than mathematical facts,
pictures, examples and surprises.

\smallskip

One is tempted to consider this reaction in the light of
Lyotard's ([L]) famous definition of the postmodern condition
as  ``incredulity toward meta--narratives'' and
Tasi\'c's remark that mathematical truth
belongs to ``the most stubborn meta--narratives of Western 
culture'' ([T], p. 176).

\smallskip

In this stubbornness lies our hope.

\bigskip

\centerline{\bf References}

\medskip

[B] H.~Bloom. {\it The Western Canon.} Riverhead Books, New York, 1994.

\smallskip

[C] G.~Cantor. {\it Beitr\"age zur Begr\"undung der transfiniten
Mengenlehre.} Math. Ann. 46 (1895), 481-- 512;
49 (1897), 207 --246.

\smallskip

[D] J.~W.~Dauben. {\it Georg Cantor: his mathematics and philosophy 
of the infinite.} Princeton UP, Princeton NJ, 1990.

\smallskip

[DeSPST] S.~Dehaene, E.~Spelke, P.~Pinet, R.~Stanescu,
S.~Tsivkin. {\it Sources of mathematical thinking:
behavioral and brain--imaging evidence.} Science, 7 May 1999, vol. 284,
970--974.

\smallskip

[F] C.~Freiling. {\it Axioms of symmetry: throwing darts
at the real line.} J.~Symb.~Logic,
51 (1986), 190--200.

\smallskip

[GaJ] M.~Garey, D.~Johnson. {\it Computers and intractability: a guide
to the theory of $NP$--completeness.} W.~H.~Freeman and Co., San--Francisco, 1979.

\smallskip

[L] J.-F.~Lyotard. {\it The postmodern condition: a report on knowledge.}
University of Minneapolis Press, Minneapolis, 1984.

\smallskip

[Man] Yu.~Manin. {\it Classical computing, quantum computing,
and Shor's factoring algorithm.} S\'eminaire Bourbaki, no. 862 (June 1999),
Ast\'erisque, vol 266, 2000, 375--404.

\smallskip

[Mas] M.~Mashaal. {\it Bourbaki.} Pour la Science, No 2, 2000.

\smallskip

[Mum] D.~Mumford. {\it The dawning of the age of stochasticity.}
In: Mathematics: Frontiers and Perspectives 2000, AMS, 1999,  197--218.

\smallskip

[PI] W.~Purkert, H.~J.~Ilgauds. {\it Georg Cantor, 1845 -- 1918.}
Birkh\"auser Verlag, Basel --  Boston -- Stuttgart, 1987.

\smallskip

[T] V. Tasi\'c. {\it Mathematics and the roots of postmodern thought.}
Oxford UP, 2001.

\bigskip

\centerline{\bf Appendix: Chronology of Cantor's life and mathematics}

\smallskip

\centerline{\it (following [PI] and [D])}

\medskip

March 3, 1845: Born in St Petersburg, Russia.

\medskip

1856: Family moves to Wiesbaden, Germany.

\medskip

1862 -- 1867: Cantor studies in Z\"urich, Berlin,
G\"ottingen and again Berlin.

\medskip

1867 -- 1869: First publications in number theory, quadratic forms.

\medskip

1869: Habilitation in the Halle University.

\medskip

1870 -- 1872: Works on convergence of trigonometric series.

\smallskip

1872 -- 1879: Existence of different
magnitudes of infinity, bijections $\bold{R}\to \bold{R}^n$,
studies of relations between continuity and dimension.

\smallskip

November 29, 1873: Cantor asked in a letter to Dedekind whether
there might exist a bijection between $\bold{N}$ and
$\bold{R}$ ([D], p.49). Shortly after Christmas
he found his diagonal procedure ([D], p. 51 etc.).

\smallskip

1874: First publication on set theory.

\smallskip

1879 -- 1884: Publication of the series
{\it \"Uber unendliche lineare Punktmannichfaltigkeiten}.

\smallskip

1883: {\it Grundlagen einer allgemeinen Mannigfaltigkeitslehre. Ein
mathematisch--philosophischer Versuch in der Lehre
des Unendlichen.}

\smallskip

May 1884: First nervous breakdown, after a successful and enjoyable trip to Paris:
depression  lasting till Fall: [D]. p. 282.

\smallskip

1884--85: Contact with Catholic theologians, encouragement from them,
but isolation in Halle. [D], p. 146: ``... early in 1885, Mittag--Leffler seemed
to close the last door on Cantor's hopes for understanding and encouragement among mathematicians.''

\smallskip

September 18, 1890: Foundation of the German Mathematical Society;
Cantor becomes its first President.

\smallskip

1891: Kronecker died.

\smallskip

1895 -- 1895: {\it ``Beitr\"age zur Begr\"undung der transfiniten 
Mengenlehre''}
(Cantor's last
major mathematical publication).

\smallskip

1897: The first ICM. Set theory is very visible.

\smallskip

1897: ``Burali--Forti [...] was the first mathematician to make public the paradoxes of 
transfinite set theory'' ([D]). He argued that {\it all} ordinals, if any pair of them
is comparable, would form an Ordinal which is greater than itself. He concluded
that not all ordinals are comparable. Cantor instead believed that all ordinals
do not form an ordinal, just as all sets do not form a set.

\smallskip

1899: Hospitalizations in Halle Nervenklinik before and after the death of son Rudolph.

\smallskip

1902--1903, winter term: Hospitalization.

\smallskip

Oct. 1907 -- June 1908: Hospitalization.

\smallskip

Sept. 1911 -- June 1912: Hospitalization.

\smallskip

1915: Celebration of the 70th anniversary of Cantor's birth, on a national level because of WW I.

\smallskip

May 1917 --  Jan. 6, 1918: Hospitalization; Cantor dies at Halle Klinik.

\vskip1cm

{\it Max--Planck--Institut f\"ur Mathematik, Bonn, Germany

and Northwestern University, Evanston, USA.}

\medskip

{\it e-mail:} manin\@mpim-bonn.mpg.de

\enddocument